\documentclass[12pt,draft]{article}

\usepackage{amsfonts}
\usepackage{amsmath}
\usepackage{amssymb}
\usepackage{amscd}
\usepackage{amsthm}
\usepackage{indentfirst}
\usepackage[hmargin=3.7cm,vmargin=3cm]{geometry}

\newtheorem{thm}{Theorem}[section]
\newtheorem{lemma}[thm]{Lemma}
\newtheorem{prop}[thm]{Proposition}
\newtheorem{coroll}[thm]{Corollary}

\theoremstyle{definition}

\newtheorem{defin}[thm]{Definition}
\newtheorem{rem}[thm]{Remark}
\newtheorem{exam}[thm]{Example}

\newtheorem*{prf}{Proof}

\newcommand{\cA}{{\mathcal{A}}}
\newcommand{\cB}{{\mathcal{B}}}

\newcommand{\cF}{{\mathcal{F}}}

\newcommand{\cM}{{\mathcal{M}}}
\newcommand{\cN}{{\mathcal{N}}}

\newcommand{\fc}{{:\ }}
\newcommand{\ve}{\varepsilon}

\newcommand{\ol}{\overline}

\newcommand{\tb}{\textbf}

\DeclareMathOperator{\id}{id}

\DeclareMathOperator{\Int}{int}

\DeclareMathOperator{\ed}{ed}

\title{Isotopy-invariant topological measures on closed orientable surfaces of higher genus}
\author{Frol Zapolsky}

\renewcommand{\labelenumi}{(\roman{enumi})}

\begin{document}

\maketitle

\begin{abstract}Given a closed orientable surface \(\Sigma\) of genus at least two, we establish an affine isomorphism between the convex compact set of isotopy-invariant topological measures on \(\Sigma\) and the convex compact set of additive functions on the set of isotopy classes of certain subsurfaces of \(\Sigma\). We then construct such additive functions, and thus isotopy-invariant topological measures, from probability measures on \(\Sigma\) together with some additional data. The map associating topological measures to probability measures is affine and continuous. Certain Dirac measures map to simple topological measures, while the topological measures due to Py and Rosenberg arise from the normalized Euler characteristic.
\end{abstract}

\section{Introduction and results}

The purpose of this paper is to give a construction of topological measures on closed orientable surfaces of genus at least two. These are invariant under isotopies, and as such are genuine, that is, non-subadditive. The motivation behind such a construction is twofold. First, to generalize the Py-Rosenberg topological measures (see below), and second, to provide simple topological measures on orientable surfaces other than the sphere and the torus. A result by the author \cite{Z_quasi-states_poisson_br_surfaces} has as a condition the existence of such a simple topological measure.

Moreover, the bulk of the constructed topological measures is on spaces of Aarnes genus \(0\) (see \cite{Aa_qm_construct} for the definition), which in the case of compact CW complexes is equivalent to vanishing first integral cohomology, and so the only orientable closed surface with this property is the sphere \(S^2\). Apart from that, there are the works of Grubb \cite{Grubb_irred_partitions_constr_qm} for general spaces of Aarnes genus \(1\) and the works of Knudsen \cite{Kn_extreme_qm,Kn_qm_torus} concerning topological measures on the torus. In any case, there has been little to no work of constructing topological measures on surfaces of higher genus, and the present paper hopefully helps to fill the gap.

The idea behind the construction is rather simple. First, a few definitions. Fix a closed orientable surface \(\Sigma\) of genus at least two. By a submanifold of \(\Sigma\) we mean a compact \(2\)-dimensional submanifold with boundary (the submanifold may be empty). A submanifold is normal if its boundary circles are non-contractible and pairwise non-homotopic. Denote the isotopy class of a normal submanfiold \(V\) by \([V]\) and the set of isotopy classes of normal submanifolds by \(\cN\). Two classes \(N,N'\in\cN\) are called disjoint if they contain disjoint representatives. There is then a well-defined sum \(N+N'\). Call a function \(\nu\fc\cN\to[0,1]\) additive if \(\nu(N+N')= \nu(N) + \nu(N')\) for disjoint \(N,N'\) and \(\nu([\Sigma])=1\). The set of additive functions is convex and compact (with respect to a suitable topology). An isotopy-invariant topological measure \(\tau\) on \(\Sigma\) gives rise to an additive function \(\nu_\tau\). The first result is

\begin{thm}\label{affine_homeo_tms_and_addv_fcns}The map \(\tau\mapsto\nu_\tau\) is an affine homeomorphism between the set of isotopy-invariant topological measures on \(\Sigma\) and the set of additive functions.
\end{thm}

What is the inverse map? For an additive function \(\nu\) the corresponding topological measure \(\tau\) is constructed as follows. A topological measure is reconstructible from its values on submanifolds. It turns out that there is a unique map \(\Phi\), which we call normalization, which associates to each submanifold \(W\subset\Sigma\) an isotopy class of normal submanifolds \(\Phi(W)\), by eliminating any disks and annuli from both \(W\) and its complement. Then \(\tau(W)=\nu(\Phi(W))\).

Next we construct additive functions. Choose a Riemannian metric of constant curvature \(-1\) on \(\Sigma\). Then in each class \(N\in\cN\) there is a unique geodesic subsurface \(U\), that is, one with geodesic boundary components; moreover, if \(N,N'\in\cN\) are disjoint, then the corresponding geodesic subsurfaces \(U,U'\) intersect only along the boundary. It now suffices to produce a function \(\mu\) which assigns a number \(\in[0,1]\) to each geodesic subsurface, such that for two subsurfaces \(U,U'\) intersecting only along the boundary it is true that \(\mu(U\cup U') = \mu(U) + \mu(U')\), and such that \(\mu(\Sigma) = 1\). An example of such a function is furnished by a probability measure which vanishes on all simple closed geodesics (which exists, since there are only countably many such geodesics).

To record the general result, let \(\cM\), \(\cF\) denote the collections of probability measures on \(\Sigma\) and of additive functions, respectively. Let \(\Delta\) denote the totality of the following data: a Riemannian metric on \(\Sigma\) of constant curvature \(-1\), and for each simple closed geodesic \(\gamma \subset \Sigma\), a coorientation \(c_\gamma\) and a weight \(w_\gamma\in[0,1]\). Then

\begin{thm}\label{measures_to_addv_fcns}To \(\Delta\) there corresponds an affine continuous map \(\cM \to \cF\).
\end{thm}

We mentioned above that simple topological measures on \(\Sigma\) are of special significance. To this end, let \(\cM_0\) denote the collection of probability measures on \(\Sigma\) which vanish on all simple closed geodesics.

\begin{coroll}The above map \(\cM\to\cF\) sends Dirac measures from \(\cM_0\) to simple topological measures.
\end{coroll}

It is interesting to determine the extent to which this map is not injective. See the discussion at the end of section \ref{tms_and_addv_fcns}.

Another example of an additive function is the normalized Euler characteristic \(\nu_0([V])=\frac{-\chi(V)}{2g-2}\), where \(g\) is the genus of \(\Sigma\). In \cite{Py_quasi_morphismes_et_inv_Calabi}, Py constructed a remarkable functional on the set of smooth functions on \(\Sigma\), which, as was shown by Rosenberg \cite{Ros_qm_construct}, is (the restriction of) a quasi-state. We call it the Py-Rosenberg quasi-state.

\begin{thm}\label{Py_Rosenberg_qs}The topological measure corresponding to \(\nu_0\) is the topological measure representing the Py-Rosen\-berg quasi-state.
\end{thm}

The paper is organized as follows. Section \ref{tms_on_mfds} deals with generalities of topological measures on closed manifolds. In section \ref{normalization} we define normal manifolds, their isotopy classes, additive functions, and prove properties thereof. We also define and establish the uniqueness of the map \(\Phi\) mentioned above. In section \ref{tms_and_addv_fcns} we prove Theorems \ref{affine_homeo_tms_and_addv_fcns}, \ref{measures_to_addv_fcns}, construct examples of additive functions, and prove Theorem \ref{Py_Rosenberg_qs}. In section \ref{exist_normalizations} we prove that \(\Phi\) exists.

\section{Topological measures on manifolds}\label{tms_on_mfds}

We denote the closure, the interior and the complement of a subset \(A\) of a topological space by \(\overline A\), \(\Int A\) and \(A^c\), respectively.

\subsection{Topological measures}

Let us remind the axioms for a (normalized) topological measure \(\tau\) on a compact space \(X\) (\cite{Aa_quasi_states_and_quasi_measures}). Denote by \(\cA\) the set of subsets \(A\subset X\) which are open or closed. Then \(\tau \fc \cA \to [0,1]\) is a topological measure if
\renewcommand{\labelenumi}{(\roman{enumi})}

\begin{enumerate}
\item \(\tau(X)=1\);
\item \(\tau(K)+\tau(K^c)=1\) for a compact \(K \in \cA\);
\item \(\tau(K)\leq\tau(K')\) for compact \(K,K'\in\cA\) with \(K\subset K'\);
\item \(\tau(K\cup K') = \tau(K)+\tau(K')\) for disjoint compact \(K,K'\in\cA\);
\item \(\tau(K)+\sup\{\tau(K')\,|\,K'\in\cA\text{ compact with }K\cap K' = \varnothing\}=1\).
\end{enumerate}

The last axiom differs from the usual one, namely
\[\tau(U)=\sup\{\tau(K)\,|\,K\in\cA\text{ compact with }K\subset U\}\]
for an open \(U\in \cA\), but together with axiom (ii) is equivalent to it.

\subsection{Reduction to manifolds and isotopies}

We are interested in topological measures on surfaces and to this end we shall make a reduction to the case of manifolds. Throughout this subsection \(M\) is a closed manifold.

\begin{defin}If \(W\subset M\) is a compact submanifold with boundary such that \(\dim W = \dim M\), we shall call \(W\) simply a submanifold (\(W\) may be empty). We denote by the collection of all such submanifolds by \(\cB\), or by \(\cB(M)\) if there is a need to stress the dependence on \(M\). The boundary of \(W\) is denoted by \(\partial W\). A connected component of \(\partial W\) is referred to as a boundary component of \(W\).
\end{defin}

There are certain operations on \(\cB\). If \(W,W'\in\cB\) are disjoint, then \(W\uplus W' \in \cB\). The inversion of \(V\in\cB\) is \(V^i:=\ol{V^c}\in\cB\); note that \(V^i \cap V = \partial V = \partial (V^i)\).

\begin{rem}By an isotopy on \(M\) we mean a smooth family \(\phi_t\) of diffeomorphisms of \(M\) with \(\phi_0=\id_M\).
\end{rem}

\begin{lemma}\label{eps_neighd_isotopy}Choose a Riemannian metric on \(M\) and let \(W \in \cB\). Then there is \(\ve_0 > 0\) such that for all \(\ve \in [0,\ve_0)\) the closed \(\ve\)-neighborhood \(W_\ve\) of \(W\) is again a submanifold, and such that for any \(\eta > 0\) there is a smooth isotopy \(\{\phi_t\}_{t\in[0,\ve]}\) of \(M\) supported in \(W_{\ve+\eta}\cap (W^i)_\eta\) satisfying \(\phi_t(W)=W_t\).
\end{lemma}

\begin{prf}Use the tubular neighborhood theorem applied to \(\partial W\). \qed \end{prf}

This allows to define two more operations on \(\cB\). For \(W\in\cB\) the closed \(\ve\)-neighbor\-hood \(W_\ve\) of \(W\) and \(W^\ve:=(W_\ve)^i\) are in \(\cB\) for all \(\ve\geq 0\) sufficiently small. Note that \(W^0 = W^i\).

\begin{rem}If \(W\in\cB\), then for any \(W'\in\cB\) disjoint from it there is \(\ve>0\) small enough such that \(W'\subset W^\ve\).
\end{rem}

Since we are on a manifold, compact sets can be approximated by submanifolds, and disjoint compact sets can be separated by disjoint submanifolds. More precisely:

\begin{lemma}\label{approx_cpt_sets_by_sbmfds}
\begin{enumerate}
\item Let \(K \subset M\) be a compact subset. Then for every open \(U \supset K\) there is a submanifold \(W\) such that \(K \subset W \subset U\);
\item Let \(K,K' \subset M\) be disjoint compact subsets. Then there are disjoint submanifolds \(W,W'\) with \(K \subset W\) and \(K' \subset W'\).
\end{enumerate}
\end{lemma}

\begin{prf}(i) Let \(F_0 \fc M \to [0,1]\) be a continuous function with \(F_0|_K \equiv 0\) and \(F_0|_{U^c} \equiv 1\). Take a smooth function \(F\) on \(M\) with \(\|F - F_0\| < \frac 1 3\) (in the uniform norm). Let \(\alpha \in (\frac 1 3, \frac 2 3)\) be a regular value of \(F\). Then \(W = \{F \leq \alpha\}\) does the job.

(ii) Since \(M\) is a normal space, there are disjoint open \(U,U' \subset M\) such that \(K\subset U\) and \(K'\subset U'\). Now use (i) to obtain submanifolds \(W,W'\) with \(K\subset W \subset U\) and \(K' \subset W' \subset U'\). \qed
\end{prf}

Hence we obtain:

\begin{coroll}\label{tm_uniquely_detd_by_sbmfds}
A topological measure \(\tau\) is uniquely determined by its restriction to \(\cB\).
\end{coroll}

\begin{prf}
It is enough to show that the values of \(\tau\) on compact subsets of \(M\) are uniquely defined by the restriction. But if \(K \subset M\) is compact, then by regularity
\[\tau(K) = \inf \{\tau(U) \,|\, U \supset K,\,U\text{ open}\}\,.\]
By Lemma \ref{approx_cpt_sets_by_sbmfds}, for any \(U \supset K\) there exists a submanifold \(W\) between \(K\) and \(U\). By monotonicity \(\tau(K) \leq \tau(W) \leq \tau(U)\), which implies that \(\tau(K) = \inf \{\tau(U) \,|\, U \supset K,\,U\text{ open}\} \geq \inf \{\tau(W) \,|\, W\in\cB\text{ with }W \supset K\}\). The reverse inequality is obvious, so we can restore the value \(\tau(K)\) as
\[\tau(K)=\inf \{\tau(W) \,|\, W\in\cB\text{ with }W \supset K\}\,.\]
\qed
\end{prf}

\begin{rem}Note the for an open subset \(O\subset M\) we have
\[\tau(O)=\sup\{\tau(W)\,|\,W\in\cB,\,W\subset O\}\,.\]
\end{rem}

Let \(\tau\) be a topological measure on \(M\). Put \(\tau_0:=\tau|_\cB\). The function \(\tau_0 \fc \cB \to [0,1]\) has the following properties:

\begin{enumerate}
\item \(\tau_0(M) = 1\);
\item \(\tau_0\) is monotone;
\item \(\tau_0\) is additive with respect to (finite) disjoint unions;
\item \(\tau_0(W) + \sup \{\tau_0(W') \,|\, W' \in \cB: W \cap W' = \varnothing\} = 1 \).
\end{enumerate}

Indeed, (i-iii) are clear, while (iv) follows from an argument similar to that of Corollary \ref{tm_uniquely_detd_by_sbmfds}. It turns out that \(\tau \mapsto \tau|_\cB\) is a bijection between the set of topological measures on \(M\) and the set of functions \(\cB \to [0,1]\) satisfying (i-iv) above:

\begin{prop}\label{reconstruction}
Suppose that we are given a function \(\tau_0 \fc \cB \to [0,1]\) satisfying the above properties. Then it is the restriction to \(\cB\) of a unique topological measure \(\tau\).
\end{prop}

\begin{prf}
Since uniqueness has already been shown, let us prove the existence.

For a compact \(K\) and an open \(U\) put
\[\tau(K) = \inf\{\tau_0(W)\,|\,W \in \cB\text{ with }W \supset K\}\quad \text{and} \quad\tau(U) = 1 - \tau(U^c)\,.\]
We have immediately: \(\tau(K) + \tau(K^c) = 1\) for any compact \(K\); \(\tau(K) \leq \tau(K')\) for compact \(K \subset K'\); \(\tau(K \uplus K') = \tau(K) + \tau(K')\) for any compact disjoint \(K,K'\). Also, clearly \(\tau(M) = 1\). It remains to show regularity. We have
\begin{align*}
\tau(K) &= \inf\{\tau_0(W)\,|\,W \in \cB\text{ with }W \supset K\}\\
&=\inf_{W\in\cB,\,W \supset K} \Big(1 - \sup_{W' \in \cB,\,W \cap W' = \varnothing} \tau_0(W')\Big)\\
&=\inf_{W\in\cB,\,W \supset K}\quad \inf_{W' \in \cB,\,W \cap W' = \varnothing} (1 - \tau_0(W'))\\
&=\inf_{W,W'\in\cB,\,W \supset K,\,W \cap W' = \varnothing} (1 - \tau_0(W'))\\
&=\inf_{W'\in\cB,W'\cap K = \varnothing} (1 - \tau_0(W'))\,,
\end{align*}
where the last equality follows from an argument similar to that of Corollary \ref{tm_uniquely_detd_by_sbmfds}. Now subtracting both sides of the final equality \(\tau(K) = \inf(1 - \tau_0(W'))\) from \(1\), we get for \(U = K^c\):
\[\tau(U) = \sup_{W' \in \cB,\,W' \subset U}\tau_0(W')\,.\]
Since \emph{a priori}
\[\tau(U) \geq \sup_{K\text{ compact with }K \subset U} \tau(K) \geq \sup_{W' \in \cB,\,W' \subset U}\tau_0(W')\,,\]
we get
\[\tau(U) = \sup_{K\text{ compact with }K \subset U} \tau(K)\,.\]
\qed
\end{prf}

We shall actually produce functions \(\tau_0\fc\cB(\Sigma)\to[0,1]\) satisfying (i-iv) above for a closed orientable surface \(\Sigma\) of genus at least two.

\subsection{Submanifolds}

In this subsection we assume that \(M\) is a closed \emph{connected} manifold.

\begin{rem}Note that \(W\in\cB\) such that \(W\neq\varnothing,M\) necessarily has nonempty boundary. Put differently, a nonempty boundaryless submanifold of \(M\) equals \(M\).
\end{rem}

\begin{defin}Let \(V,W\in\cB\). We say that \(V,W\) intersect well if every connected component of \(\partial V\) is either disjoint from \(\partial W\) or coincides with a connected component of \(\partial W\). In this case every connected component of \(\partial V\cap \partial W\) is a connected component of both \(\partial V\) and \(\partial W\). We say that \(V,W\) are compatible if they intersect well and also \(\ol{\Int V \cap \Int W} \supset \partial V \cap \partial W\). We say that \(V\) dominates \(W\) is they are compatible and \(\partial V \supset \partial W\).
\end{defin}

\begin{rem}
Let \(V,W\in\cB\). If \(V,W\) intersect well, then \(V\cup W \in \cB\). If they are compatible, then \(\ol{V\triangle W}, \ol{V-W}\in\cB\), and \(\partial\ol{V\triangle W} = \partial V \triangle \partial W\). If in addition \(V\supset W\) then \(\partial\ol{V-W}=\partial V \triangle \partial W\). If \(V\) dominates \(W\) then \(\partial\ol{V\triangle W} = \partial V - \partial W\).
\end{rem}

\begin{lemma}Let \(V,W\in\cB\) and let \(V\) dominate \(W\). Let \(\gamma\) be a boundary component of \(W\). Then the connected component of \(W\) containing \(\gamma\) contains the connected component of \(V\) containing \(\gamma\).
\end{lemma}

\begin{prf}Denote the connected components of \(V,W\) containing \(\gamma\) by \(S,T\), respectively. We need to show that \(T\supset S\). The boundary of \(T\), \(\partial T\), separates \(M\), \(\Int T\) being one of the connected components of the complement \(M-\partial T\). By domination, \(\Int S\) is disjoint from \(\partial T\), and because it is connected, it must be contained in one of the connected components of \(M-\partial T\). Again by domination, \(\Int S \cap \Int T \neq \varnothing\) because \(\ol{\Int S \cap \Int T} \supset \gamma\). Thus \(\Int S \subset \Int T\) and the claim follows. \qed
\end{prf}

\begin{coroll}\label{mutual_domination_implies_equality}Let \(V,W\in\cB\) dominate each other. If \(\partial V\neq\varnothing\), then \(V=W\).
\end{coroll}

\begin{prf}Let \(V=\biguplus_iS_i\) and \(W=\biguplus_jT_j\) be the decompositions into connected components. Since \(M\) is connected and \(\partial V = \partial W\) is nonempty, each one of the submanifolds \(S_i\) and \(T_j\) has boundary. Now if \(\gamma\) is a boundary component of \(V\) and \(S_i\), \(T_j\) are the respective components containing \(\gamma\), then by the lemma, \(S_i=T_j\). Define \(V'=V-S_i\) and \(W'=W-T_j\). These submanifolds again dominate each other and have nonempty boundary. By induction on the number of connected components in \(V\) we have \(V'=W'\) and so \(V=W\). \qed
\end{prf}

\begin{defin}For a closed codimension \(1\) submanifold \(H\subset M\) a coorientation is an orientation of the conormal bundle \(N^*_MH\).
\end{defin}

\begin{rem}If \(W \in \cB\), then the boundary \(\partial W\) acquires a canonical coorientation, namely, the outward one. If \(W'\) is another submanifold, then \(W\) dominates \(W'\) if and only if \(\partial W' \subset \partial W\) and the coorientations of \(\partial W\) and \(\partial W'\) agree. It follows from Corollary \ref{mutual_domination_implies_equality} that a nonempty submanifold is completely determined by its cooriented boundary. For future use, let us also note that the isotopy class of a (nonempty) submanifold is uniquely determined by the isotopy class of its cooriented boundary.
\end{rem}

\section{Normalization}\label{normalization}

For the rest of the paper \(\Sigma\) stands for a closed orientable surface of genus \(g \geq 2\).

\subsection{Topology of surfaces}

\begin{defin}
By a curve \(\gamma \subset \Sigma\) we mean a smoothly embedded circle; \(\gamma\) is contractible if there is a parametrization \(\ol\gamma\fc S^1\to\Sigma\) of \(\gamma\) such that \(\ol\gamma\) is null-homotopic; a curve \(\gamma'\) is homotopic to \(\gamma\) if there is a parametrization \(\ol{\gamma'}\) of \(\gamma'\) such that \(\ol\gamma\) and \(\ol{\gamma'}\) are freely homotopic. A disk in \(\Sigma\) is a smoothly embedded closed disk. An annulus in \(\Sigma\) is a smoothly embedded annulus \(S^1\times[0,1]\), where we assume that the boundary components are non-contractible.

\end{defin}

Note that a boundary component of submanifold of \(\Sigma\) is a curve.

The following are standard facts from surface topology, see for example \cite{Bu_geom_spectr_cpt_riem_surf}.

\begin{lemma}
\begin{enumerate}
\item A contractible curve \(\gamma \subset \Sigma\) is the boundary of a unique disk, denoted by \(D(\gamma)\);
\item Two disjoint non-contractible homotopic curves \(\gamma,\gamma' \subset \Sigma\) form the boundary of a unique annulus, denoted by \(A(\gamma,\gamma')\).
\item If \(\gamma_1,\dots,\gamma_r\) and \(\gamma'_1,\dots,\gamma'_r\) are two collections of disjoint non-con\-tract\-ible pairwise non-homotopic curves such that \(\gamma_i,\gamma'_i\) are homotopic for all \(i\), then there is an isotopy \(\phi\) such that \(\gamma'_i=\phi(\gamma)\) for all \(i\).
\qed
\end{enumerate}
\end{lemma}

This allows to introduce a partial order on the set of contractible curves: \(\gamma\leq\gamma'\) if \(D(\gamma)\subset D(\gamma')\). In particular, since the set of contractible boundary components of a submanifold is finite, there is a finite number of components which are maximal with respect to this partial order.

Similarly, there is a partial order on the set of unordered pairs of disjoint non-contractible homotopic curves, defined by \([\gamma,\gamma']\leq [\delta,\delta']\) if \(A(\gamma,\gamma')\subset A(\delta,\delta')\) (\([a,b]\) denotes the unordered pair consisting of elements \(a,b\)). Given three non-contractible homotopic curves \(\gamma,\gamma',\gamma''\) it makes sense to say that \(\gamma\) is between \(\gamma'\) and \(\gamma''\), meaning that \(\gamma \subset \Int A(\gamma',\gamma'')\). More generally, for a finite collection of non-contractible homotopic curves \(\gamma_1,\dots,\gamma_k\) it makes sense to say that they are listed in undirected order if \(\gamma_j\) is between \(\gamma_{j-1}\) and \(\gamma_{j+1}\) for all \(j=2,\dots,k-1\).


There is a finite number of maximal pairs of boundary components of a submanifold.

\subsection{Normal submanifolds}

\begin{defin}A submanifold \(V\in\cB\) is called normal if its boundary components are pairwise non-homotopic.
\end{defin}

\begin{rem}A submanifold \(V\) is normal if and only every component of both it and its inversion admit a decomposition into pairs-of-pants, see \cite{Bu_geom_spectr_cpt_riem_surf}. The inversion of a normal submanifold is again normal.
\end{rem}

\begin{defin}Let \(V,V'\in\cB\) be normal and disjoint. There is a unique normal \(V''\) associated to \(V \uplus V'\), see subsection \ref{eliminating_annuli}. Call this \(V''\) the sum of \(V\) and \(V'\): \(V''=V+V'\).
\end{defin}

\begin{lemma}\label{normal_sbmfd_subtraction}Let \(V,V'\in\cB\) be normal with \(V\subset V'\). Then either \(V\) is isotopic to \(V'\) or there is an isotopy \(\psi\) and a nonempty normal \(V'' \in \cB\) unique up to isotopy such that \(V''\) is disjoint from \(\psi(V)\) and \(V'=V''+\psi(V)\).
\end{lemma}

\begin{prf} If \(V\) is isotopic to \(V'\), we are done. Otherwise, by applying an isotopy to \(V\) if necessary, and denoting the new submanifold still by \(V\), we may arrange that every boundary component \(\gamma\) of \(V\) is a boundary component of \(V'\) unless \(\gamma\) is not homotopic to any boundary component of \(V'\), in which case we may assume that \(\gamma\) is disjoint from \(\partial V'\); call such a boundary component of \(V\) interior to \(V'\).

Clearly, the manifold \(V'':=\ol{V'-V}\) is normal, and applying an additional isotopy to \(V\) we may displace a little every boundary component of \(V\) interior to \(V'\) such that it is disjoint from \(V''\). It follows that \(V'=V''+V\). This shows that such \(V''\) exists. The isotopy class of the cooriented boundary of \(V''\) is that of \((\partial V'-\partial V)\uplus(\partial V - \partial V')\), where the first part is taken with the coorientation of \(\partial V'\) while the second part is given the coorientation opposite to that of \(\partial V\). This shows that the isotopy class of \(V''\), being uniquely determined by the isotopy class of its cooriented boundary (since \(V''\) is nonempty), is also uniquely determined by \(V,V'\).\qed
\end{prf}

We use the standard notation \(\chi(X)\) for the Euler characteristic of a topological space \(X\).

\begin{defin}For a normal submanifold \(V\) define its complexity \(n(V)\) to be \(-\chi(V)\).
\end{defin}

\begin{rem}Euler characteristic satisfies the so-called inclusion-exclusion principle, which states that for two subsets \(A,B\) of a topological space it is true, under some conditions, that \(\chi(A\cup B)=\chi(A)+\chi(B)-\chi(A\cap B)\). These conditions are certainly satisfied if \(A,B\) are two submanifolds of \(\Sigma\) which intersect well and only along the boundary. In this case the intersection is a finite disjoint union of circles, and so has zero Euler characteristic. Since the Euler characteristic of a pair-of-pants is \(-1\), it follows that the complexity of a normal submanifold \(V\) equals the number of pairs-of-pants in a decomposition of \(V\). In particular, \(n(V)=0\) if and only if \(V=\varnothing\).
\end{rem}

\begin{lemma}
\begin{enumerate}
\item If \(V,V'\) are normal and disjoint then \(n(V+V')=n(V)+n(V')\);
\item If \(V,V'\) are normal with \(V\subset V'\) then \(n(V)\leq n(V')\) with equality if and only if \(V\) is isotopic to \(V'\).
\end{enumerate}
\end{lemma}

\begin{prf} (i) Follows from the remark above.

(ii) If \(V\) is isotopic to \(V'\), we are done. Otherwise it follows from Lemma \ref{normal_sbmfd_subtraction} that \(V'=V''+\psi(V)\), so
\[n(V')=n(V'')+n(\psi(V))=n(V'')+n(V)>n(V)\,,\]
since \(V''\) is nonempty and hence has nonzero complexity.\qed
\end{prf}

\subsection{Isotopy classes of normal submanifolds and additive functions}

\begin{defin}For a normal submanifold \(V\in\cB\) let \([V]\) denote its isotopy class. Let \(\cN\) denote the set of isotopy classes of normal submanifolds. If \(N,N'\in\cN\), we say that \(N,N'\) are disjoint if there are normal disjoint \(V,V'\in\cB\) such that \(N=[V],N'=[V']\); in this case we define the sum \(N+N':=[V+V']\). We say that \(N\) is contained in \(N'\), and write \(N\subset N'\), if there are normal \(V,V'\in\cB\) such that \(N=[V],N'=[V']\) and \(V\subset V'\). The class \([\varnothing]\) is called the empty class and is denoted by \(\varnothing\). The class \([\Sigma]\) will be denoted \(\Sigma\). Lastly, for \(N=[V]\) put \(N^i:=[V^i]\).
\end{defin}

\begin{rem}Note that \(N+\varnothing = N\) for any \(N\in\cN\), and that \(N,N^i\) are always disjoint with \(N+N^i=\Sigma\).
\end{rem}

\begin{lemma}The sum of disjoint isotopy classes is well-defined.
\end{lemma}

\begin{prf} If both \(N,N'\) are the empty class, there is nothing to prove. Otherwise we only need to show that the isotopy class of the cooriented boundary of the sum \(V+V'\) is independent of the choice of disjoint representatives \(V,V'\) of \(N,N'\). The boundary of \(V+V'\) is \(\partial V \cup \partial V' - \bigcup_i(\gamma_i\cup \gamma'_i)\), where \(\gamma_i\subset \partial V,\gamma'_i\subset \partial V'\) are the pairs of homotopic boundary components. We see that the isotopy class of this cooriented boundary is uniquely defined by the classes of \(V,V'\), and hence that the class of the sum \(V+V'\) is independent of the choice of \(V,V'\) within their isotopy classes.\qed
\end{prf}

\begin{lemma}Let \(N,N'\in\cN\) be such that \(N\subset N'\). Then there is a unique \(N''\in\cN\) disjoint from \(N\) such that \(N'=N+N''\).
\end{lemma}

\begin{prf} Choose representatives \(V,V'\) of \(N,N'\) such that \(V\subset V'\). Then Lemma \ref{normal_sbmfd_subtraction} yields an isotopy \(\psi\) and a normal submanifold \(V''\) unique up to isotopy, such that \(V''\) is disjoint from \(\psi(V)\) and \(V=V''+\psi(V')\). It follows that \(N''=[V'']\) is disjoint from \(N\) and \(N'=N''+N\), and that the class of \(V''\) is uniquely determined by \(N,N'\).\qed
\end{prf}

\begin{defin}The complexity of a class \([V]\in\cN\) is \(n([V]):=n(V)\).
\end{defin}

\begin{lemma}
\begin{enumerate}
\item If \(N,N'\in\cN\) are disjoint, then \(n(N+N')=n(N)+n(N')\);
\item If \(N,N'\in\cN\) are such that \(N\subset N'\), then \(n(N)\leq n(N')\) with equality if and only if \(N=N'\).
\end{enumerate}
\end{lemma}

\begin{prf} (i) Follows from the corresponding property of the complexity of normal submanifolds.

(ii) Let \(N''\) be the unique class such that \(N'=N''+N\). Then \(n(N')=n(N'')+n(N)\geq n(N)\). The equality is attained if and only if \(n(N'')=0\) which is equivalent to \(N''=\varnothing\) in which case \(N'=N+\varnothing=N\). \qed
\end{prf}

\begin{coroll}The containment relation on \(\cN\) is a partial order.
\end{coroll}

\begin{prf}The only nontrivial property is symmetry. So let \(N,N'\in\cN\) be such that \(N\subset N'\) and \(N'\subset N\). Then \(\nu(N)=\nu(N')\) and so \(N=N'\). \qed
\end{prf}

\begin{defin}An additive function is a function \(\nu\fc\cN\to[0,1]\) such that
\begin{enumerate}
\item If \(N,N'\in\cN\) are disjoint, then \(\nu(N+N')=\nu(N)+\nu(N')\);
\item \(\nu(\Sigma)=1\).
\end{enumerate}
An additive function is called simple if it only takes values \(0\) and \(1\).
\end{defin}

\begin{rem}Note that an additive function is monotone. Indeed, if \(N\subset N'\), let \(N''\) be such that \(N'=N+N''\); then \(\nu(N')=\nu(N)+\nu(N'')\geq\nu(N)\).
\end{rem}

A convex combination of additive functions is again such, so the set \(\cF\) of additive functions is affine. Topologize \(\cF\) as follows: a net \(\nu_i\in\cF\) converges to \(\nu\in\cF\) if and only if
\[\liminf_i \nu_i(N) \geq \nu(N)\]
for any \(N\in\cN\). Since \(\nu_i\) and \(\nu\) are additive and normalized, we have
\begin{multline*}1-\nu(N)=\nu(N^i)\leq\liminf_i\nu_i(N^i)\leq\limsup_i\nu_i(N^i)=\\=1-\liminf_i \nu(N)\leq 1 - \nu(N)\,,\end{multline*}
so the topology thus defined is that of pointwise convergence. The evaluation map \(\cF\to[0,1]^\cN,\,\nu\mapsto(\nu(N))_N\) is then an embedding with closed image and so \(\cF\) is compact. The subspace of simple additive functions is closed, and so is a compact space on its own.

\subsection{Uniqueness of normalization}

\begin{defin}An isotopy-invariant map \(\Phi\fc\cB\to\cN\) is called a normalization if for \(W,W'\in\cB\):
\begin{enumerate}
\item \(W\cap W'=\varnothing\) implies that \(\Phi(W),\Phi(W')\) are disjoint and \(\Phi(W\uplus W')=\Phi(W)+\Phi(W')\);
\item \(W\subset W'\) implies \(\Phi(W)\subset \Phi(W')\);
\item \(\Phi(W^i)=\Phi(W)^i\);
\item if \(W\) is normal then \(\Phi(W)=[W]\).
\end{enumerate}
\end{defin}

We shall eventually prove that there is a unique normalization. In this subsection we prove the uniqueness.

\begin{lemma}\label{vanishing_on_disks_and_annuli}If \(W\in\cB\) is a disk or an annulus then \(\Phi(W)=\varnothing\) for any normalization \(\Phi\), and \(\tau(W)=0\) for any isotopy-invariant topological measure \(\tau\) on \(\Sigma\).
\end{lemma}

\begin{prf}Since any disk (or annulus) can be isotoped off itself, and the disjoint union of two disks (two isotopic annuli) can be isotoped into a third disk (annulus), we have that the class \(\Phi(W)\) is disjoint from itself and moreover \(\Phi(W)+\Phi(W)\subset\Phi(W)\). Hence \(2 n(\Phi(W))\leq n(\Phi(W))\) which implies \(n(\Phi(W))=0\), and this can happen if and only if \(\Phi(W)= \varnothing\). \qed
\end{prf}

\begin{thm}\label{uniqueness_normalizations}A normalization \(\Phi\), if it exists, is unique.
\end{thm}

\begin{prf}Let \(W\in\cB\). If \(W\) is normal, \(\Phi(W)=[W]\), and so is uniquely determined. Otherwise there are two cases. Case 1: \(W\) has a contractible boundary component. Let \(\gamma\) be minimal such a component. If \(\gamma\) is exterior to \(W\) (see subsection \ref{eliminating_disks}), then for \(W'=W-D(\gamma)\) we have \(W=W'\uplus D(\gamma)\). By Lemma \ref{vanishing_on_disks_and_annuli}, \(\Phi(D(\gamma))=\varnothing\), and so \(\Phi(W)=\Phi(W')\uplus\Phi(D(\gamma))=\Phi(W')\). But \(W'\) has fewer boundary components and so by induction \(\Phi(W)=\Phi(W')\) is uniquely determined. If \(\gamma\) is interior to \(W\), then it is exterior to \(W^i\) and hence \(\Phi(W)=\Phi(W^i)^i\) is again uniquely determined. Case 2: \(W\) has no contractible boundary components, by there are two homotopic boundary components. Let \([\gamma,\gamma']\) be minimal such a pair. Then it is either inner or outer with respect to \(W\) (see subsection \ref{eliminating_annuli}). If it is outer, for \(W'=W-A(\gamma,\gamma')\) we have \(W=W'\uplus A(\gamma,\gamma')\) and by the argument of Case 1 \(\Phi(W)\) is uniquely determined; if the pair is inner, consider again \(W^i\), with respect to which it is outer and reduce to the previous case. \qed
\end{prf}
\section{Topological measures and additive functions}\label{tms_and_addv_fcns}
\subsection{From topological measures to additive functions}

If \(\tau\) is an isotopy-invariant topological measure, then we can define \(\nu\fc\cN\to[0,1]\) by \(\nu([V]):=\tau(V)\).

\begin{prop}This \(\nu\) is an additive function.
\end{prop}

\begin{prf} We have \(\nu(\Sigma)=\tau(\Sigma)=1\). Now let \(N,N'\) be disjoint. Then there are normal disjoint \(V,V'\) with \(N=[V],N'=[V']\), and \(V+V'=V\uplus V'\uplus A\), where \(A\) is a finite union of open annuli with disjoint closures. Lemma \ref{vanishing_on_disks_and_annuli} implies that \(\tau(A)=0\), therefore
\begin{multline*}\nu(N+N')=\nu([V+V'])=\tau(V+V')=\tau(V)+\tau(V')+\tau(A)=\\=\nu([V])+\nu([V'])=\nu(N)+\nu(N')\,.\end{multline*}
\qed
\end{prf}

\begin{prop}The map \(\tau\mapsto\nu\) is affine and continuous and transfers simple topological measures to simple additive functions.
\end{prop}

\begin{prf} Affinity is obvious, so let us prove continuity. Let \(\tau_i\to\tau\) be a convergent net and let \(\nu_i,\nu\) be the corresponding additive functions. Let \(N\in\cN\), and choose \(V\in N\). Since an isotopy-invariant topological measure vanishes on embedded circles, \(\tau_i(V)=\tau_i(O)\) for all \(i\) and also \(\tau(V)=\tau(O)\) where \(O=\Int V\). Then
\[\liminf \nu_i(N)=\liminf_i\tau_i(O)\geq \tau(O) = \nu(N)\,,\]
as required. The last assertion is clear.\qed
\end{prf}

\subsection{From additive functions to topological measures}

If \(\nu\) is an additive function, put \(\tau_0(W):=\nu(\Phi(W))\) for \(W \in\cB\), where \(\Phi\) is the normalization map.

\begin{prop}This \(\tau_0\) extends to a unique isotopy-invariant topological measure \(\tau\) and the additive function corresponding to \(\tau\) is \(\nu\).
\end{prop}

\begin{prf} By Proposition \ref{reconstruction} we must prove that \(\tau_0\) is (i) normalized, (ii) additive, (iii) monotone and (iv) regular.

(i) \(\tau_0(\Sigma)=\nu(\Sigma)=1\).

(ii) If \(W,W'\) are disjoint, then so are \(\Phi(W)\) and \(\Phi(W')\) are and \(\Phi(W\uplus W')=\Phi(W)+\Phi(W')\), hence
\begin{multline*}\tau_0(W\uplus W')=\nu(\Phi(W\uplus W'))=\nu(\Phi(W)+\Phi(W'))=\\=\nu(\Phi(W))+\nu(\Phi(W')) = \tau_0(W)+\tau_0(W')\,.\end{multline*}

(iii) If \(W\subset W'\), then \(\Phi(W)\subset\Phi(W')\), hence
\[\tau_0(W)=\nu(\Phi(W)) \leq \nu(\Phi(W')) = \tau_0(W')\,,\]
because additive functions are monotone.

(iv) We need to show that \(\tau_0(W) + \sup \{\tau_0(W') \,|\, W' \in \cB: W \cap W' = \varnothing\} = 1\). For sufficiently small \(\ve_0>0\) every \(W^\ve\) with \(\ve\in(0,\ve_0)\) is a submanifold disjoint from \(W\), thus
\[\sup \{\tau_0(W') \,|\, W' \in \cB: W \cap W' = \varnothing\} \geq \sup_{\ve \in (0,\ve_0)} \tau_0(W^\ve)\,.\]
Since \(W\) is a compact submanifold, any \(W'\in \cB\) disjoint from it is contained in \(W^\ve\) for some \(\ve\in(0,\ve_0)\), which implies an inequality reverse to that above. Thus we need to show \(\tau_0(W) + \sup_{\ve \in (0,\ve_0)} \tau_0(W^\ve) = 1\). But for \(\ve \geq 0\) sufficiently small \(W^\ve\) is isotopic to \(W^i\) and so \(\Phi(W^\ve)=\Phi(W^i)=\Phi(W)^i\), whence \(\sup_{\ve \in (0,\ve_0)} \tau_0(W^\ve)=\nu(\Phi(W)^i)\). The classes \(\Phi(W)\) and \(\Phi(W)^i\) are disjoint and \(\Phi(W)+\Phi(W)^i=\Sigma\), thus
\[\tau_0(W) + \sup_{\ve \in (0,\ve_0)} \tau_0(W^\ve) = \nu(\Phi(W))+\nu(\Phi(W)^i)=\nu(\Sigma)=1\,.\]
The last assertion, as well as the invariance of \(\tau\) under isotopies, are obvious.
\qed
\end{prf}

Next we prove that the above map \(\nu\mapsto\tau\) is continuous, and for this we require two lemmas.

\begin{lemma}Let \(O\subset \Sigma\) be an open subset. Then the family \(\cB_O:=\{W\in\cB\,|\, W \subset O\}\), viewed as a partially ordered set with respect to inclusion, is directed, that is for any \(W,W'\in\cB_O\) there is \(W''\in\cB_O\) such that \(W\subset W''\) and \(W'\subset W''\).
\end{lemma}

\begin{prf} Let \(K=W\cup W'\). Then \(K\) is a compact subset contained in \(O\) and so there is a submanifold \(W''\) with \(K\subset W''\subset O\).\qed
\end{prf}

\begin{lemma}\label{maximal_sbmfd_open_set_add_fcn}Let \(O\subset \Sigma\) be an open subset. There is \(W_0\in \cB_O\) such that for \(W\in\cB_O\) we have \(\Phi(W_0)\supset\Phi(W)\).
\end{lemma}

\begin{prf}The complexity is a bounded integer-valued function on \(\cN\), and so there is \(W_0\in\cB_O\) such that \(n(\Phi(W_0))\geq n(\Phi(W))\) for any \(W\in\cB_O\). We claim that \(W_0\) satisfies the desired property. Indeed, let \(W\in\cB_O\) and let \(W''\in\cB_O\) be such that \(W_0\cup W \subset W''\). Then \(\Phi(W'')\supset\Phi(W_0)\), hence by monotonicity of the complexity \(n(\Phi(W''))\geq n(\Phi(W_0))\), but \(n(\Phi(W_0))\geq n(\Phi(W''))\) by the choice of \(W_0\) and so \(\Phi(W_0)=\Phi(W'')\supset \Phi(W)\). \qed
\end{prf}

\begin{prop}The map \(\nu\mapsto\tau\) is affine and continuous and maps simple additive functions to simple topological measures.
\end{prop}

\begin{prf} Affinity being obvious, let us prove continuity. Let \(\nu_i\to\nu\) be a convergent net of additive functions, and let \(\tau_i,\tau\) be the corresponding topological measures. We need to show that for any open \(O\subset \Sigma\) we have
\[\liminf_i\tau_i(O)\geq \tau(O)\,.\]
Let \(W_0\in\cB_O\) be as in Lemma \ref{maximal_sbmfd_open_set_add_fcn}. Then for any additive function \(\pi\) we have \(\sup\{\pi(\Phi(W))\,|\,W\in\cB_O\}=\pi(\Phi(W_0))\). This follows from the fact that the left-hand side is bounded from above by \(\pi(\Phi(W_0))\), since additive functions are monotone, while the reverse inequality is obvious. Let \(\sigma\) be the topological measure corresponding to \(\pi\). Then we have \(\sigma(O)=\pi(\Phi(W_0))\). Indeed, by regularity
\[\sigma(O)=\sup\{\sigma(W)\,|\,W\in\cB_O\}=\sup\{\pi(\Phi(W))\,|\,W\in\cB_O\}=\pi(\Phi(W_0))\,.\]
It follows that
\[\liminf_i\tau_i(O)=\liminf_i\nu_i(\Phi(W_0))\geq \nu(\Phi(W_0))=\tau(O)\,.\]
The last claim is obvious.
\qed
\end{prf}

\subsection{Examples of additive functions}

\begin{exam}The normalized complexity \(\nu:=\dfrac {n(\cdot)} {2g-2}\) is an additive function (\(g\) is the genus of \(\Sigma\)). Let \(\tau\) denote the corresponding topological measure. Then

\begin{thm}The topological measure \(\tau\) corresponds to the Py-Rosenberg quasi-state \(\zeta\) on \(\Sigma\).
\end{thm}

\begin{prf}Denote the quasi-state corresponding to \(\tau\) by \(\zeta_\tau\). By the Lipschitz continuity of a quasi-state it suffices to show that \(\zeta\) and \(\zeta_\tau\) agree on the set of generic Morse functions on \(\Sigma\).

Recall the definition of \(\zeta\). Let \(F\in C^\infty(M)\) be a generic Morse function. There is a notion of an essential critical point of \(F\). There are \(2g-2\) of them. Let \(\alpha_1<\dots<\alpha_{2g-2}\) be the critical values of \(F\) corresponding to essential critical points. Then
\[\zeta(F)=\frac 1 {2g-2}\sum_{i=1}^{2g-2}\alpha_i\,.\]

From the work of Rosenberg \cite{Ros_qm_construct} it follows that \(\alpha\) is the critical value of \(F\) corresponding to an essential critical point if and only if for all \(\ve>0\) small enough \(\Phi({\{F\leq \alpha + \ve\}})\) differs from \(\Phi({\{F\leq \alpha - \ve\}})\) by one pair-of-pants. The function \(b_F\), defined by \(b_F(t):=\tau(\{F\leq t\})\) then equals \(\frac i {2g-2}\) for \(t\in[\alpha_i,\alpha_{i+1})\) for \(i=1,\dots,2g-3\), \(0\) for \(t < \alpha_1\) and \(1\) for \(t\geq \alpha_{2g-2}\). By definition
\begin{align*}\zeta_\tau(F) &= \max F - \int_{\min F}^{\max F}b_F(t)\,dt\\
&= \alpha_1+\frac 1 {2g-2}\sum_{i=1}^{2g-3}(2g-2-i)(\alpha_{i+1}-\alpha_i)\\
&= \frac 1 {2g-2}\sum_{i=1}^{2g-2}\alpha_i=\zeta(F)\,.
\end{align*}
\qed
\end{prf}

\end{exam}

To construct more general additive functions, we need to measure classes of normal submanifolds. A universal method is to choose canonical (in some sense) representatives of those classes and then to measure them. To do this, choose a Riemannian metric of constant curvature \(-1\) on \(\Sigma\). Call a curve on \(\Sigma\) geodesic if it has a geodesic parametrization, and call a submanifold geodesic if each one of its boundary components is a geodesic curve. It is known that in any free homotopy class of curves on \(\Sigma\) there is a unique geodesic. More generally, any finite collection of non-contractible pairwise non-homotopic curves is isotopic to a unique geodesic such collection, hence any normal submanifold is isotopic to a unique geodesic submanifold. Thus we have a way of choosing representatives from classes in \(\cN\), which we call geodesification (for the lack of a better name); the geodesification of a class \(N\in\cN\) is denoted by \(\widetilde N\).

\begin{defin}If \(U,U'\) are two geodesic submanifolds, we say that they are essentially disjoint if they intersect only along the boundary: \(U\cap U' \subset \partial U \cap \partial U'\).
\end{defin}

\begin{rem}If geodesic \(U,U'\) are essentially disjoint, then their intersection consists of a finite number of geodesic circles, each one being a component of both \(\partial U\) and \(\partial U'\). Therefore they intersect well and their union \(U\cup U'\) is then also geodesic. The inversion of a geodesic submanifold is again such.
\end{rem}

Geodesification has the following properties:

\begin{prop}Let \(N,N'\in\cN\). Then
\begin{enumerate}
\item If \(N,N'\) are disjoint, it implies that \(\widetilde N\) and \(\widetilde {N'}\) are essentially disjoint; moreover, \(\widetilde{N+N'}=\widetilde N \cup \widetilde{N'}\);
\item \(N\subset N'\Rightarrow \widetilde N \subset \widetilde{N'}\);
\item \(\widetilde{N^i}=\big(\widetilde N\big)^i\).
\end{enumerate}
\end{prop}

\begin{prf}(i) If both \(N,N'\) are empty, there is nothing to prove. Otherwise, let \(V,V'\) be disjoint representatives of \(N,N'\). There is an isotopy \(\psi\) fixing the boundary components of \(V\) which are not homotopic to any of the boundary components of \(V'\), such that \(\psi(V)\cap V'\) is precisely the union of all boundary components \(\gamma'\) of \(V'\) such that there is a boundary component \(\gamma\) of \(V\) which is homotopic to \(\gamma'\). It follows that \(V+V'=\psi(V)\cup V'\), since a nonempty submanifold is uniquely determined by its cooriented boundary, which in this case is \(\partial \psi(V) \triangle \partial V\). 

Let \(C=\partial \psi(V) \cup \partial V'\). Then \(C\) is isotopic to a unique collection of pairwise disjoint geodesic curves by an isotopy \(\phi\). Then both \(\phi(\psi(V))\) and \(\phi(V')\) are geodesic and essentially disjoint. By the uniqueness of geodesification \(\phi(\psi(V))=\widetilde N,\phi(V')=\widetilde{N'}\). On the other hand, \(\phi(V+V')=\phi(\psi(V)\cup\phi(V'))=\phi(\psi(V))\cup\psi(V')=\widetilde N\cup \widetilde {N'}\) is geodesic and is isotopic to \(V+V'\), so \(\widetilde{N+N'}=\phi(V+V')=\widetilde N\cup \widetilde {N'}\).

(ii) Let \(N''\) be the unique class disjoint from \(N\) such that \(N'=N+N''\). Then \(\widetilde{N'}=\widetilde{N}\cup\widetilde{N''}\supset\widetilde N\).

(iii) Clear.
\qed
\end{prf}

\begin{exam}\label{constr_addv_fcns_via_prob_meas_0}
There are countably many free homotopy classes of curves on \(\Sigma\), so there are at most countably many geodesic curves. 
In particular, there are probability measures on \(\Sigma\) which vanish on all geodesic curves. So choose such a probability measure \(\mu\) and put \(\nu(N):=\mu(\widetilde N)\).
\end{exam}

\begin{lemma}\(\nu\) is an additive function.
\end{lemma}

\begin{prf}(i) \(\nu(\Sigma)=\mu(\Sigma)=1\).

(ii) Let \(N,N'\in\cN\) be disjoint. Then \(\widetilde N\) and \(\widetilde {N'}\) are essentially disjoint, that is \(\widetilde N\cap\widetilde {N'}\) is a finite union of geodesic curves and so has \(\mu\)-measure zero. Then
\[\nu(N+N')=\mu(\widetilde N\cup\widetilde {N'})=\mu(\widetilde N)+\mu(\widetilde {N'})=\nu(N)+\nu(N')\,.\]
\qed
\end{prf}

\begin{prop}The map \(\mu\mapsto\nu\) is affine, continuous, and maps Dirac measures to simple additive functions.
\end{prop}

\begin{prf}The first and the last assertions are clear. To prove continuity, let \(N\in\cN\). We need to show that if \(\mu_i\) is a net of probability measures as above converging to \(\mu\), then the corresponding net of additive functions \(\nu_i\) converges to \(\nu\). For any compact subset \(K\subset\Sigma\) we have
\[\limsup_i \mu_i(K)\leq \mu(K)\,,\]
and so
\begin{multline*}\liminf_i\nu_i(N)=1-\limsup_i\nu(N^i) =1-\limsup_i\mu_i((\widetilde N)^i)\geq\\
\geq 1-\mu((\widetilde N)^i) =1-(1-\mu(\widetilde N))=\nu(N)\,.
\end{multline*}
\qed
\end{prf}

\begin{exam}\label{constr_addv_fcns_via_prob_meas}This example generalizes the previous one. Let \(\mu\) be arbitrary probability measure on \(\Sigma\). If \(\gamma\) is a curve, it has two possible coorientations. For a cooriented curve \(\delta\) let \(c(\delta)\) denote the coorientation. A boundary component of a submanifold of \(\Sigma\), has a natural coorientation, namely, the outward one. Now, for each geodesic \(\gamma\), choose a coorientation \(c_\gamma\) and a weight \(w_\gamma\in[0,1]\). If \(\gamma',\gamma''\) are two cooriented versions of the same curve, define the product \(c(\gamma')c(\gamma'')\) to be \(1\) if the coorientations agree and \(-1\) otherwise. For \(w \in [0,1]\) and \(\sigma = \pm 1\) put \(\langle w,\sigma\rangle\) to be \(w\) if \(\sigma = 1\) and \(1-w\) if \(\sigma = -1\). Now if \(N\in\cN\) and \(U:=\widetilde N\) has the cooriented boundary components \(\gamma_1,\dots,\gamma_r\), put
\[\nu(N):=\mu(\Int U)+\sum_{i=1}^r \langle w_{\gamma_i},c(\gamma_i)c_{\gamma_i}\rangle\mu(\gamma_i)\,.\]
It is a matter of routine verification that this defines an additive function on \(\cN\). The definition is tailored so that when two geodesic submanifolds are essentially disjoint, the contribution of the intersection is exactly its \(\mu\)-measure.
\end{exam}

\begin{rem}In \cite{BS_geodesics_on_surfaces} it is shown that the union of all\footnote{That is, simple and closed, according to our terminology.} geodesic curves on \(\Sigma\) is nowhere dense and has Hausdorff dimension \(1\). In particular, its complement \(C\), which is exactly the intersection of the interiors of all pairs-of-pants in \(\Sigma\), contains large open connected subsets. Each connected component of \(C\) gives rise to exactly one additive function as in Example \ref{constr_addv_fcns_via_prob_meas_0}, that is, Dirac measures with supports inside the same connected component of \(C\) map to the same additive function. It is an interesting problem to determine, for instance, whether the set of connected components of \(C\) is uncountable. Another way of getting simple additive functions is to take a probability measure supported on geodesic curves and to assign to each geodesic a weight which is either \(0\) or \(1\), as in Example\ref{constr_addv_fcns_via_prob_meas}. How many simple topological measures are obtained in this way is unknown.
\end{rem}

\section{Proof of existence of normalizations}\label{exist_normalizations}
The idea of proof is actually contained in the proof of uniqueness, Theorem \ref{uniqueness_normalizations}. Arguing by induction on the number of boundary components of a submanifold \(W\), we add or remove a disk for each contractible boundary component of \(W\), thus reducing their number; for each pair of non-contractible homotopic boundary components of \(W\) we analogously add or remove an annulus, again reducing the number of boundary components. At each step we have to keep track of disjoint unions and containment, as well as inversions. In this section we carry out the details of the construction.

Therefore, there are two steps in this construction. The first one, called eliminating disks, disposes of any contractible boundary components of a given submanifold \(W\in\cB\). There is a way to define the resulting submanifold uniquely in terms of \(W\). Eliminating disks commutes with isotopies, preserves disjoint unions and containment, and also taking inversions.

The second step is called eliminating annuli. This is a procedure which associates a normal submanifold to a given one with no contractible boundary components. There is no natural way to define the result of eliminating annuli uniquely so that it satisfy all the desired properties, but to any submanifold we associate a finite number of possibilities, in particular, there are unique minimal and maximal associated normal submanifolds. All the normal submanifolds arising from the given one are isotopic. If two given submanifolds are disjoint, we can associate to them disjoint normal submanifolds, and if they contain one another, again, the normal counterparts may be chosen so that the containment is preserved. The same is true of inversion.

The normalization \(\Phi(W)\) of a given submanifold \(W\) is then defined as the (uniquely determined!) isotopy class of a normal submanifold associated with to the result of eliminating disks in \(W\). The desired properties of \(\Phi\) then follow from those of eliminating disks and annuli.

\subsection{Eliminating disks}\label{eliminating_disks}

\begin{defin}An element \(U \in \cB\) is called disk-free, if \(\partial U\) has no contractible components.
\end{defin}

This is equivalent to requiring that no component of \(U\) or of \(U^i\) be contained in a closed disk.

\begin{rem}A submanifold is contained in a (finite) disjoint union of closed disks if and only if each one of its boundary components is contractible. If two such submanifolds intersect well, their union is again contained in a disjoint union of closed disks. Note that such a submanifold is necessarily distinct from the whole surface \(\Sigma\).
\end{rem}

\begin{prop}Given \(W \in \cB\), there is a unique disk-free submanifold \(U\), such that
\begin{enumerate}
\item \(W\) dominates \(U\);
\item \(\ol{U \triangle W}\) is contained in a finite union of disjoint closed disks.
\end{enumerate}
\end{prop}

First let us give a definition, which also will be needed later.
\begin{defin}If \(\gamma \subset \partial W\) is a contractible connected component, then if \(W\) dominates \(D(\gamma)\), we say that \(\gamma\) is exterior to \(W\); otherwise we say that \(\gamma\) is interior to \(W\).
\end{defin}

\begin{prf}First, let us check uniqueness. Let \(U,V\) satisfy the required conditions. It follows that \(U,V\) dominate each other and so if they have nonempty boundary, they are equal by Corollary \ref{mutual_domination_implies_equality}. If their boundaries are empty, we only have to rule out the case \(U=\varnothing\), \(V=\Sigma\). Assume that this is the case. Then \(W=\ol{\varnothing \triangle W}=\ol{U\triangle W}\) is contained in a disjoint union of disks. On the other hand, \(W^i=\ol{\Sigma - M}=\ol{V\triangle W}\) is also contained in a disjoint union of disks. Since \(W,W^i\) intersect well, their union \(W\cup W^i = \Sigma\) again is contained in a disjoint union of disks, which is impossible.

Let us now construct such a submanifold out of \(W\). Let \(\gamma_1,\dots,\gamma_r\) be the maximal contractible components of \(\partial W\), the first \(p\) being interior to \(W\) and the remaining \(r-p\) being exterior to \(W\). Then the submanifold
\[U:=\Big(W\cup\bigcup_{j=1}^p D(\gamma_j)\Big)-\bigcup_{j=p+1}^r D(\gamma_j)\]
satisfies the required conditions. \qed

\end{prf}

The manifold \(U\) thus obtained is said to be the result of eliminating disks in \(W\) and is denoted by \(\ed W\).

\begin{rem}Note that a disk-free \(U\) equals \(\ed W\) if and only if \(W\) dominates \(U\) and every component of \(\partial W - \partial U\) is contractible.
\end{rem}

In order to establish the necessary properties of eliminating disks, we need the
\begin{lemma}Let \(\gamma\) be a contractible curve and let \(W\) be a submanifold contained in \(D(\gamma)\). Then every maximal boundary component of \(W\) is exterior to \(W\).
\end{lemma}

\begin{prf}If \(W\) is empty, there is nothing to prove. Otherwise let \(\delta_1\) be a maximal boundary component of \(W\) and let \(W_0\) be the connected component of \(W\) containing \(\delta\). Let \(\delta_1,\delta_2,\dots,\delta_r\) be the boundary components of \(W_0\). If \(\delta_1\) is interior to \(W\), let \(W'_0:=W_0\cup D(\delta_1)\). If \(\delta_j\subset\partial W'_0\) is exterior to \(W'_0\) for some \(j=2,\dots,r\), then by connectedness \(W_0\subset W'_0\subset D(\delta_j)\), hence \(\delta_1\subset W_0\subset D(\delta_j)\), contradicting the maximality of \(\delta_1\). Hence all \(\delta_j\) are interior to \(W'_0\). It follows that \(W'_0\cup\bigcup_{j=2}^r D(\delta_j)\) is a nonempty boundaryless submanifold of \(\Sigma\) contained in \(D(\gamma)\), which is absurd. This contradiction shows that \(\delta_1\) is exterior to \(W\).\qed
\end{prf}

\begin{prop}Let \(V,W \in \cB\). We have the following properties of eliminating disks:
\begin{enumerate}
\item if \(V,W\in\cB\) are disjoint, then so are \(\ed V, \ed W\), and \(\ed(V\uplus W)=\ed V \uplus \ed W\);
\item if \(V\subset W\), then \(\ed V \subset \ed W\);
\item \(\ed(W^i) = (\ed W)^i\).
\end{enumerate}
\end{prop}

\begin{prf}(i) By symmetry, the only thing we need to show is that if \(\gamma\) is a maximal contractible component of \(\partial V\) which is interior to \(V\), then \(D(\gamma)\cap \ed W = \varnothing\). Note that if \(\delta \subset \partial W\) is a contractible component, then, since \(V,W\) are disjoint, either \(D(\gamma),D(\delta)\) are disjoint, \(\delta \subset \Int D(\gamma)\), or \(\gamma\subset\Int D(\delta)\). Let us rule out the last case. Consider the connected component \(V_0\) of \(V\) containing \(\gamma\). By connectedness, it lies in \(\Int D(\delta)\), and since \(\gamma\) is interior to \(V_0\), by the above lemma it contradicts the maximality of \(\gamma\).

Now let \(\delta\subset\partial W\) be a maximal contractible component. The relevant case is \(\delta \subset \Int D(\gamma)\). By the lemma \(\delta\) is exterior, and so is eliminated from \(W\) in the eliminating disks algorithm.

It remains to prove that \(\ed(V\uplus W)=\ed V \uplus \ed W\). Since \(\ed V \uplus \ed W\) is disk-free, \(V\uplus W\) dominates \(\ed V \uplus \ed W\) and \(\partial (V \uplus W) - \partial (\ed V \uplus \ed W) = (\partial V - \partial (\ed V)) \uplus (\partial W - \partial (\ed W))\) has only contractible components, by the above remark, \(\ed V \uplus \ed W\) must be equal \(\ed(V\uplus W)\).

(ii) We need to show that if \(\gamma\subset\partial V\) is a maximal connected component interior to \(V\), then \(D(\gamma)\subset \ed W\), or equivalently, \(D(\gamma)\cap \ed (W^i) = \varnothing\). Now \(\ed (W^i)\) can meet \(V\) only at boundary points. Let \(\delta \subset \partial W\) be a maximal contractible component. There are three possible cases: (a) \(D(\delta)\cap D(\gamma) = \varnothing\); (b) \(\delta \subset D(\gamma)\) and (c) \(\gamma \subset D(\delta)\). Case (a) does not concern us; case (c) is impossible, because we would get that the connected component of \(V\) containing \(\gamma\) lies in \(D(\delta)\), contradicting either the maximality of \(\gamma\) or the fact that it is interior to \(V\). So we are left with case (b). Since \(W^i\) meets \(V\) only along the boundary, it follows that \(\delta\) is exterior to \(W^i\): just look at the submanifold \(W^i \cap D(\gamma)\) and use the above lemma. Hence \(D(\gamma)\) does not meet \(\ed W^i\), because \(D(\delta)\) is eliminated from \(W^i\) in the algorithm.

(iii) Since \((\ed W)^i\) is disk-free, \(W^i\) dominates it and every component of \(\partial (W^i) - \partial (\ed W)^i = \partial W - \partial (\ed W)\) is contractible, by the remark, \((\ed W)^i = \ed (W^i)\).\qed
\end{prf}

\subsection{Eliminating annuli}\label{eliminating_annuli}

\begin{defin} We say that a normal submanifold \(V\) is associated to a disk-free submanifold \(W\) if (i) \(W\) dominates \(V\), and (ii) \(\ol{V \triangle W}\) is a disjoint union of annuli.
\end{defin}

\begin{rem}If we are given two submanifolds both of which are finite disjoint unions of closed annuli, then their union is again such, provided they intersect well. Note that such a submanifold is necessarily distinct from the whole surface \(\Sigma\).
\end{rem}

Associated normal submanifolds are not unique, but there are relatively few of them for any given submanifold, as can be seen from the following

\begin{prop}\label{assoc_normal_uniq_detd_by_bdry}Let \(W\) be disk-free and let \(U,V\) be normal submanifolds associated to \(W\) such that \(\partial U = \partial V\). Then \(U=V\).
\end{prop}

\begin{prf}\(U,V\) dominate each other, so unless \(\partial U=\partial V=\varnothing\), \(U=V\). In case \(U,V\) have no boundary, we only have to rule out the case \(U=\varnothing,V=\Sigma\). But then \(W=\ol{\varnothing\triangle W}=\ol{U\triangle W}\) is a disjoint union of annuli, as well as \(W^i=\ol{\Sigma\triangle W}=\ol{V\triangle W}\). Since \(W,W^i\) intersect well it follows that \(\Sigma=W\cup W^i\) is a disjoint union of annuli, which is impossible. \qed
\end{prf}

In particular, there are finitely many associated normal submanifolds for any given disk-free \(W\in\cB\). We shall see that any two of them are isotopic, and in fact, we shall get a complete description of all of them.

Let us establish the structure of disk-free submanifolds. First, a definition.
\begin{defin}Let \(W\in \cB\) and let \(\gamma,\gamma'\) be two disjoint non-contractible homotopic boundary components of \(W\). We call the unordered pair \([\gamma,\gamma']\) inner, intermediate, or outer with respect to \(W\) if \(\ol{\Int W \cap A(\gamma,\gamma')}\) contains none of \(\gamma,\gamma'\), one of them, or both of them, respectively. Of the curves \(\gamma,\gamma'\) call those contained in \(\ol{\Int W \cap A(\gamma,\gamma')}\) outer and those which are not contained in this set inner.\footnote{Note that a curve being inner or outer depends on the pair to which it belongs!} If \(\delta\) is the only connected component of \(\partial W\) in its homotopy class, \(\delta\) is called intermediate.
\end{defin}

Now let \(W \in \cB\) be disk-free. The boundary \(\partial W\) (if non-empty) falls into connected components
\[\partial W = \biguplus_{j=1}^r\biguplus_{k=1}^{k_j}\gamma_{j,k}\,,\]
where \(\gamma_{j,k}\) and \(\gamma_{j,k'}\) are homotopic for all \(j,k,k'\), and no two \(\gamma_{j,k}\), \(\gamma_{j',k'}\) are homotopic for \(j\neq j'\). It follows that for each \(j\) the curves \(\{\gamma_{j,k}\}_k\) can be arranged in undirected order, and without loss of generality \(\gamma_{j,1},\dots,\gamma_{j,k_j}\) is the order. If \(k_j \geq 2\), then the pair \([\gamma_{j,1},\gamma_{j,k_j}]\) is maximal.

If the pair \([\gamma_{j,1},\gamma_{j,k_j}]\) is outer, then \(k_j\) is even, and the connected components of \(W\) inside \(A(\gamma_{j,1},\gamma_{j,k_j})\) are precisely \(A(\gamma_{j,1},\gamma_{j,2}),\dots, A(\gamma_{j,k_j-1},\gamma_{j,k_j})\). If \([\gamma_{j,1},\gamma_{j,k_j}]\) is inner, then \(k_j\) is even, and the connected components of \(W\) inside \(A(\gamma_{j,1},\gamma_{j,k_j})\) are precisely \(A(\gamma_{j,2},\gamma_{j,3}),\dots, A(\gamma_{j,k_j-2},\gamma_{j,k_j-1})\). If \([\gamma_{j,1},\gamma_{j,k_j}]\) is intermediate, then without loss of generality \(\gamma_{j,1}\) is inner and \(\gamma_{j,k_j}\) is outer. It follows that \(k_j\) is odd, and so \(k_j \geq 3\), and the connected components of \(W\) inside \(A(\gamma_{j,1},\gamma_{j,k_j})\) are precisely \(A(\gamma_{j,2},\gamma_{j,3}),\dots, A(\gamma_{j,k_j-1},\gamma_{j,k_j})\).

Let us call the index \(j\) outer, inner or intermediate according to whether the pair \([\gamma_{j,1},\gamma_{j,k_j}]\) is outer, inner or intermediate, respectively. If \(k_j=1\), we shall also call \(j\) intermediate.

We shall now describe the possible normal submanifolds associated to \(W\). Assume without loss of generality that \(j=1,\dots,s\) are intermediate, that \(j=s+1,\dots, s+m\) are inner and \(j=s+m+1,\dots,r\) are outer.

Put
\[W':=\Big(W \cup \bigcup_{j=s+1}^{s+m}A(\gamma_{j,1},\gamma_{j,k_j})\Big)- \bigcup_{j=s+m+1}^r A(\gamma_{j,1},\gamma_{j,k_j})\,.\]
It is clear that
\[\partial W' = \biguplus_{j=1}^s\biguplus_{k=1}^{k_j}\gamma_{j,k}\,,\]
and so all \(j=1,\dots,s\) are intermediate (with respect to \(W'\)). Assume, again without loss of generality, that \(k_j\geq 3\) for \(j=1,\dots,p\) and \(k_j=1\) for \(j=p+1,\dots,s\).

Choose a sequence \(\underline i = (i_1,\dots,i_p)\) of odd integers such that \(1 \leq i_j \leq k_j\) for \(j=1,\dots,p\). Now put \(W'_0:=W'\) and for \(j=1,\dots,p\) define inductively \(W'_j\) to be \(W'_{j-1}\) together with all the connected components of \(W'\) lying between \(\gamma_{j,1}\) and \(\gamma_{j,i_j}\), if any, while discarding all the connected components of \(W'\) lying between \(\gamma_{j,i_j}\) and \(\gamma_{j,k_j}\), if any. It is clear that \(V_{\underline i}:=W'_p\) is a normal submanifold associated to \(W\). Actually, any normal submanifold associated to \(W\) is of this form:

\begin{prop}With notations as above, if \(V\) is a normal submanifold associated to \(W\), then there is \(\underline i\) such that \(V=V_{\underline i}\).
\end{prop}

\begin{prf}Since \(\ol{V\triangle W}\) is a submanifold whose components are all annuli, there is an even number of boundary components in each homotopy class of curves. Because \(\partial\ol{V\triangle W}\subset \partial W\), the possible homotopy classes are enumerated by the index \(j\). By virtue of the relation \(\partial\ol{V\triangle W} = \partial W - \partial V\) the parity of the number of boundary components of \(V\) in the class determined by \(j\) equals that of \(k_j\). Since \(V\) is normal, it has at most one boundary component in each homotopy class, and so the number of boundary components in the class of \(j\) equals \(0\) if \(k_j\) is even, and \(1\) is \(k_j\) is odd.

Now \(W\) dominates \(V\) and so, by the above consideration, there is a bijection between the boundary components of \(V\) and intermediate indices \(j\). For each intermediate \(j=1,\dots,p\), let \(i_j\) be the unique index such that the boundary component of \(V\) in the class determined by \(j\) is \(\gamma_{j,i_j}\). Then \(\partial V=\partial V_{\underline i}\), and by Proposition \ref{assoc_normal_uniq_detd_by_bdry}, \(V=V_{\underline i}\). \qed
\end{prf}

It is also clear that if \(\underline i = (1,\dots,1)\), then \(V_{\underline i}\) is the minimal normal submanifold associated to \(W\), that is, it is contained in any other such submanifold; and if \(\underline i = (k_1,\dots,k_p)\), then \(V_{\underline i}\) is the maximal normal submanifold associated to \(W\). From the construction it follows that for any two \(\underline i\) and \(\underline i'\) the submanifolds \(V_{\underline i}\) and  \(V_{\underline i'}\) are isotopic, because their cooriented boundaries are. Hence we have the following
\begin{prop}Let \(W\in\cB\) be disk-free and let \(V,V'\) be two normal submanifolds associated to \(W\). Then \(V,V'\) are isotopic. \qed
\end{prop}

Next we show the necessary properties of associated normal submanifolds. First, a lemma.

\begin{lemma}Let \(\gamma,\gamma'\) be two simple closed non-contractible curves and let \(W\) be a disk-free submanifold contained in \(A(\gamma,\gamma')\). Then there is \(n \geq 0\) such that the boundary components of \(W\) are \(\delta_1,\dots,\delta_{2n}\), they are all pairwise homotopic and homotopic to \(\gamma,\gamma'\), the list \(\gamma,\delta_1,\dots,\delta_{2n},\gamma'\) is in undirected order, and \(W=\biguplus_{j=1}^nA(\delta_{2j-1},\delta_{2j})\). \qed
\end{lemma}

\begin{prop}Let \(W,W'\) be disk-free. Then
\begin{enumerate}
\item If\/ \(W\cap W' = \varnothing\), then there are disjoint normal submanifolds \(V,V'\) associated to \(W,W\), respectively;
\item If\/ \(W\subset W'\), then there are normal submanifolds \(V,V'\) associated to \(W,W\), respectively, such that \(V \subset V'\);
\item Inversion sends normal submanifolds associated to \(W\) into those of \(W^i\).
\end{enumerate}
\end{prop}

\begin{prf}(i) We shall show that the minimal normal submanifolds \(V,V'\) associated to \(W,W'\) are disjoint. Since we are dealing with minimal associated normal submanifolds, intermediate maximal pairs of boundary components of \(W,W'\) are of no concern to us. Hence, by symmetry, it suffices to show that if \([\alpha,\beta]\) is a maximal boundary pair of \(W\) which is inner with respect to \(W\), then \(A(\alpha,\beta)\cap V' = \varnothing\). Assume that \(\gamma\) is a boundary component of \(W'\) homotopic to \(\alpha,\beta\). Then we claim that \(\gamma \subset \Int A(\alpha,\beta)\). Indeed, otherwise without loss of generality \(\alpha\) is between \(\gamma\) and \(\beta\). Let \(W_0\) be the connected component of \(W\) containing \(\alpha\). Since \(W\) is disjoint from \(W'\), \(W_0\) is disjoint from \(\gamma\) and so by connectedness it is true that \(W_0\subset A(\gamma,\alpha)\). But then, by the lemma, \(\partial W_0 = \alpha \uplus \alpha'\), \(\alpha'\) being necessarily between \(\gamma\) and \(\alpha\), which then contradicts the maximality of \([\alpha,\beta]\), because we obtained \([\alpha',\beta]\geq[\alpha,\beta]\).

Similarly, it follows that any connected component of \(W'\) which has a boundary component homotopic to \(\alpha,\beta\), must lie in \(\Int A(\alpha,\beta)\). Thus the maximal pair \([\gamma,\delta]\) of boundary components of \(W'\) homotopic to \(\alpha,\beta\), if it exists, satisfies \(A(\gamma,\delta)\subset\Int A(\alpha,\beta)\) and so it is outer with respect to \(W'\). It then follows that any connected components of \(W'\) lying in \(\Int A(\alpha,\beta)\) are removed in the process of eliminating annuli as described above and so \(A(\alpha,\beta)\cap V'=\varnothing\), as required.

(ii) Let \(V\) be the minimal normal submanifold associated to \(W\) and \(V'\) be the maximal normal submanifold associated to \(W'\). We claim that \(V \subset V'\). Again, intermediate maximal pairs are irrelevant and so we only have to show that if \([\alpha,\beta]\) is a maximal boundary pair of \(W\) which is inner with respect to \(W\), then \(V'\supset A(\alpha,\beta)\). By an argument similar to that of (i), any boundary component \(\gamma\) of \(W'\) is contained in \(A(\alpha,\beta)\). Let \([\gamma,\delta]\) be the maximal boundary pair of \(W'\) in the homotopy class of \(\alpha,\beta\). By maximality the pair \([\gamma,\delta]\) must be inner with respect to \(W'\). Without loss of generality \(\gamma\) lies between \(\alpha\) and \(\delta\). Then the connected component of \(W'\) containing \(\gamma\) also contains the connected component of \(W\) containing \(\alpha\) and the connected component of \(W'\) containing \(\delta\) contains the connected component of \(W\) containing \(\beta\). It follows that in the process of forming the maximal normal submanifold \(V'\) associated to \(W'\) the annulus \(A(\gamma,\delta)\) is added to \(W'\) and so \(V'\) contains the annulus \(A(\alpha,\beta)\).

(iii) Clear.\qed
\end{prf}

Finally, we would like to describe normal submanifolds associated to disjoint unions. Let \(W,W'\in\cB\) be disk-free and disjoint and let \(W'':=W\uplus W'\). Let \(V,V'\) be normal submanifolds associated to \(W,W'\), respectively, and assume that \(V,V'\) are disjoint. Let \(\partial V = \gamma_1\uplus\dots\uplus\gamma_m\) and \(\partial V' = \gamma'_1\uplus\dots\uplus\gamma'_{m'}\) be the decompositions into boundary components. Assume without loss of generality that \(\gamma_i,\gamma'_i\) are homotopic for \(i\leq p\) and that there are no other pairs of homotopic curves among the boundary components of \(V,V'\). Then for each \(i\leq p\) the boundary pair \([\gamma_i,\gamma'_i]\) of \(V\uplus V'\) is maximal, and moreover it is inner with respect to \(V\uplus V'\). From boundary considerations it follows that
\[V''=V\cup V'\cup \bigcup_{i\leq p} A(\gamma_i,\gamma'_i)\]
is the unique normal submanifold associated to \(V\uplus V'\). It is dominated by \(W\uplus W'\) and \(\partial\ol{(W\uplus W')\triangle V''}\) consists of pairs of non-contractible homotopic curves. Hence \(V''\) is a normal submanifold associated to \(W''\).

\begin{thm}Normalizations exists.
\end{thm}

\begin{prf}For \(W\in\cB\) let \(\Phi(W)\) be the class of a normal submanifold associated to \(\ed W\). It is well-defined since any two normal submanifolds associated to a disk-free submanifold are isotopic. The map \(\Phi\) is isotopy-invariant, since eliminating disks and annuli are.

(i) If \(W,W'\) are disjoint, then \(\ed(W\uplus W')=\ed W \uplus \ed W'\), and there are disjoint normal \(V,V'\) associated to \(\ed W,\ed W'\), respectively, and moreover, \(V+V'\) is normal and associated to \(\ed W \uplus \ed W'\). Hence the classes \(\Phi(W),\Phi(W')\) are disjoint and \(\Phi(W\uplus W')=[V+V']=[V]+[V']=\Phi(W)+\Phi(W')\).

(ii) If \(W\subset W'\), then \(\ed W\subset \ed W'\), and there are normal \(V,V'\) associated to \(\ed W,\ed W'\), respectively, such that \(V\subset V'\). Hence \(\Phi(W)=[V]\subset[V']=\Phi(W')\).

(iii) Clear.

(iv) The only normal submanifold associated with \(\ed V = V\), where \(V\) is normal, is \(V\) itself, so \(\Phi(V)=[V]\). \qed
\end{prf}

\noindent\textsc{Frol Zapolsky\\
School of Mathematical Sciences\\
Tel Aviv University\\
Tel Aviv 69978, Israel}\\
\texttt{zapolsky@post.tau.ac.il}


\begin{thebibliography}{Aa2}

\bibitem[Aa1]{Aa_quasi_states_and_quasi_measures} Aarnes, J. F., \textit{Quasi-states and quasi-measures}, Adv. Math. \tb{86} (1991), no. 1, 41--67.

%
\bibitem[Aa2]{Aa_qm_construct} Aarnes, J. F., \textit{Construction of non-subadditive measures and discretization of Borel measures}, Fund. Math. \tb{147} (1993), 213--237.

\bibitem[BS]{BS_geodesics_on_surfaces}Birman, J.S., Series, C., \textit{Geodesics with bounded intersection number on surfaces are sparsely distributed}, Topology \tb{24} (1985), 217--225.

\bibitem[B]{Bu_geom_spectr_cpt_riem_surf} Buser, P., \textit{Geometry and spectra of compact Riemann surfaces}, Progress in Mathematics, vol. 106, Birkh\"auser Boston, 1992.



\bibitem[G]{Grubb_irred_partitions_constr_qm} Grubb, D. J., \textit{Irreducible partitions and the construction of quasi-measures}, Trans. Amer. Math. Soc. \tb{353} (2001),  no. 5, 2059--2072 (electronic).

\bibitem[K1]{Kn_extreme_qm} Knudsen, F. F., \textit{Topology and the construction of extreme quasi-measures}, Adv. Math. \tb{120} (1996), no. 2, 302--321.

\bibitem[K2]{Kn_qm_torus} Knudsen, F. F., \textit{New topological measures on the torus}, Fund. Math. \tb{185} (2005), no. 3, 287--293.

\bibitem[P]{Py_quasi_morphismes_et_inv_Calabi} Py, P., \textit{Quasi-morphismes et invariant de Calabi},  Ann. Sci. \'Ecole Norm. Sup. (4) \tb{39}  (2006),  no. 1, 177--195.

\bibitem[R]{Ros_qm_construct} Rosenberg, M., \textit{Py-Calabi quasi-morphisms and quasi-states on orientable surfaces of higher genus}, M. Sc. thesis, arXiv: 0706.0028.


\bibitem[Z]{Z_quasi-states_poisson_br_surfaces} Zapolsky, F., \textit{Quasi-states and the Poisson bracket on surfaces}, J. Mod. Dyn.  \tb{1} (2007), no. 3, 465--475.

\end{thebibliography}
\end{document}